\documentclass[a4paper,11pt]{article}
\usepackage{amsmath}
\usepackage{amssymb}
\usepackage{amsthm}
\usepackage{mathtools}
\usepackage{suffix}
\usepackage{stmaryrd}
\usepackage{bbold}
\usepackage{breqn}
\usepackage{txfonts}
\usepackage{graphicx}
\usepackage{xfrac}
\usepackage{booktabs}

\usepackage[dvipsnames]{xcolor}
\usepackage{soul}
\usepackage{comment}
\specialcomment{coloredcomment}{\begingroup\color{blue!50}}{\endgroup}

\usepackage{accents}

\usepackage[english]{babel}
\usepackage{blindtext}

\usepackage{bm}
\bmdefine{\ba}{a}
\bmdefine{\bb}{b}
\bmdefine{\bc}{c}
\bmdefine{\bd}{d}
\bmdefine{\be}{e}
\bmdefine{\bf}{f}
\bmdefine{\bg}{g}
\bmdefine{\bj}{j}
\bmdefine{\bk}{k}
\bmdefine{\bn}{n}
\bmdefine{\bp}{p}
\bmdefine{\bq}{q}
\bmdefine{\br}{r}
\bmdefine{\bx}{x}
\bmdefine{\by}{y}
\bmdefine{\bu}{u}
\bmdefine{\bv}{v}
\bmdefine{\bw}{w}
\bmdefine{\bz}{z}
\bmdefine{\bE}{E}
\bmdefine{\bB}{B}
\bmdefine{\bC}{C}
\bmdefine{\bD}{D}
\bmdefine{\bJ}{J}
\bmdefine{\bR}{R}
\bmdefine{\bS}{S}
\bmdefine{\bT}{T}
\bmdefine{\bvarphi}{\varphi}
\bmdefine{\bphi}{\phi}
\bmdefine{\bpsi}{\psi}
\bmdefine{\bchi}{\chi}
\bmdefine{\bPsi}{\Psi}
\bmdefine{\blambda}{\lambda}
\bmdefine{\bmu}{\mu}
\bmdefine{\bnu}{\nu}

\def\xfoo#1^#2\relax#3\valign{%
\mathbf{#1}\ifx\valign#2\valign\else^{\mathbf{#2}}\fi}

\usepackage[mathscr]{euscript}

\makeatletter
\newcommand*{\defeq}{\mathrel{\rlap{%
					 \raisebox{0.3ex}{$\m@th\cdot$}}%
					 \raisebox{-0.3ex}{$\m@th\cdot$}}%
					 =}
\makeatother

\makeatletter
\newcommand*{\eqdef}{=\mathrel{\rlap{%
					 \raisebox{0.3ex}{$\m@th\cdot$}}%
					 \raisebox{-0.3ex}{$\m@th\cdot$}}%
					 }
\makeatother

\def\XXint#1#2#3{{\setbox0=\hbox{$#1{#2#3}{\int}$}
     \vcenter{\hbox{$#2#3$}}\kern-.5\wd0}}

\newcommand\wrapped[1]%
  {%
   \begin{array}{@{}l@{}}#1\end{array}%
  }
  
\theoremstyle{definition}

\theoremstyle{plain}

\theoremstyle{remark}

\theoremstyle{definition}

\DeclarePairedDelimiterX\braket[2]{\langle}{\rangle}{#1\,\delimsize\vert\,\mathopen{}#2}

\DeclarePairedDelimiterX\MeijerM[3]{\lparen}{\rparen}%
{#3\, \delimsize\vert\,\begin{smallmatrix}#1 \\ #2\end{smallmatrix}}

\newcommand\MeijerG[8][]{%
  \mathbf{G}^{\,#2,#3}_{\,#4,\,#5}\MeijerM[#1]{#6}{#7}{#8}}

\WithSuffix\newcommand\MeijerG*[7]{%
  \mathbf{G}^{\,#1,#2}_{\,#3,\,#4}\MeijerM*{#5}{#6}{#7}}

\DeclarePairedDelimiterX\FoxM[3]{\lbrack}{\rbrack}%
{#3\, \delimsize\vert\begin{smallmatrix}#1 \\ #2\end{smallmatrix}}

\newcommand\FoxH[8][]{%
  \mathbf{H}^{\,#2,#3}_{\,#4,\,#5}\FoxM[#1]{#6}{#7}{#8}}

\WithSuffix\newcommand\FoxH*[7]{%
  \mathbf{H}^{\,#1,#2}_{\,#3,\,#4}\FoxM*{#5}{#6}{#7}}

\usepackage[]{hyperref}
\hypersetup{
    bookmarks=true,
    bookmarksnumbered=true,     
    bookmarksopen=true,         
    colorlinks=true,            
    pdfstartview=Fit,           
    pdfpagemode=UseOutlines,    
    pdfpagelayout=TwoPageRight
}

\usepackage{cleveref}
\crefformat{pluralequation}{#2eqs.~(#1)#3}
\Crefformat{pluralequation}{#2Eqs.~(#1)#3}

\usepackage{authblk}

\providecommand{\keywords}[1]
{
  \small	
  \textbf{\textit{Keywords---}} #1
}

\makeatletter
\newcommand{\proofstep}[1]{%
  \par
  \addvspace{\medskipamount}%
  \noindent\textit{#1\@addpunct{.}}\enspace\ignorespaces
}
\makeatother

\newcommand{\RE}{\mathop{\mathrm{Re}}}

\newcommand{\ee}{\operatorname{e}}

\begin{document}
	
\title{Umbral methods, function factorisation and generalisation of the Fourier transform method}
\author[1]{Giuseppe Dattoli\footnote{E-mail: pinodattoli@libero.it (Giuseppe Dattoli)}}
\author[1]{Roberto Ricci\footnote{E-mail: roberto.ricci@enea.it (Roberto Ricci)}}
\author[2]{Tommaso Severati\footnote{E-mail: tom.severati@gmail.com (Tommaso Severati)}}
\affil[1]{ENEA, Nuclear Department, Frascati Research Center, Via E. Fermi 45, 00044 Frascati (Rome), Italy}
\affil[2]{University of Queensland Dept. of Mathematics and Physics Brisbane, QLD 4072, Australia}
\date{}
\setcounter{Maxaffil}{0}
\renewcommand\Affilfont{\itshape\small}
\maketitle

\begin{abstract}
	We propose a systematic way to construct trigonometric-like functions beyond the classical sine--cosine pair by factorising rational umbral operators. The guiding idea is simple: the usual trigonometric functions may be viewed as cyclic components arising from a finite factorisation, and the same principle can be extended to an $n$-fold decomposition of rational umbral expressions. For each integer $n\geq 2$, the construction produces $n$ functions which play the role of higher-order trigonometric components: their sum reconstructs the corresponding umbral function, while the individual components isolate the different cyclic sectors of its expansion.

	The construction is developed first in the formal umbral setting. The quadratic case $n=2$ gives the Gaussian trigonometric functions, in which the cosine-like component is a Gaussian and the sine-like component is its natural umbral companion. The cubic case $n=3$ yields a three-component cyclic system and shows how the same idea extends beyond the usual even--odd decomposition. These examples suggest that trigonometric factorisation is not restricted to ordinary rotations, but belongs to a broader cyclic principle in umbral calculus.

	We then reinterpret the same formal identities through the recently developed analytic umbral framework. In this second step, the cyclic components are realised by Mellin--Barnes pairings, and the root-of-unity decomposition is related to the splitting of the corresponding spectral kernel. This analytic formulation provides contour representations, local expansions, and sectorial asymptotics for the functions obtained formally. Finally, we indicate how the same cyclic kernels act on Fourier transforms. The resulting framework presents higher-order umbral trigonometric functions as natural cyclic components of factorised rational or exponential umbral operators.
\end{abstract}

\vskip0.2cm

\section{Introduction}\label{sec:introduction}

The theory of trigonometric functions beyond the elementary circular case has appeared in many different forms throughout the development of mathematics over the last three centuries.

A number of classical references, recalled in \cite{Lindqvist1995}, provide a broad historical perspective on this subject. They enter deeply into the structure of trigonometric functions and provide the technical tools needed to understand how the elementary circular functions are connected with hypergeometric and elliptic functions, and how sine- and cosine-like functions may emerge from higher transcendents, including functions belonging to the Bessel family \cite{Ferrari1981}.

Both well-known and less familiar mathematicians contributed to this development. Some of these contributions became part of standard university-level mathematical education, while others remained confined to a more specialised literature.

Euler introduced the so-called quasi-trigonometric functions \cite{Euler1795} in connection with the summation of series extending beyond the Basel problem, namely the evaluation of the infinite sum of the reciprocal squares of the positive integers. Further comments on this point will be given in the final section. Legendre and Jacobi \cite{Lawden1989} laid the foundations of the theory of elliptic functions, which arose from the inversion of elliptic integrals. Other authors, including Ferrari and Lindqvist \cite{Lindqvist1995,Ferrari1981}, discussed trigonometric functions associated with elliptic curves.

In more recent times, the emergence of new mathematical languages has renewed interest in this field (see e.g. \cite{SenturkUnal2021} and \cite{Babusci2012}). In particular, indicial umbral calculus \cite{Dattoli2017} provides a useful framework for developing trigonometric-like functions from a different point of view. The motivations behind this approach are manifold and will be discussed below.

In Ref. \cite{Dattoli2023}, Di Palma, Licciardi, and one of the present authors (G. D.) introduced Gaussian trigonometry, using the umbral theory of Gaussian functions as a benchmark. A different umbral representations of the same functions has been recently proposed by another of the present authors (R.R.) in \cite{Ricci2026}.

We briefly recall the relevant formalism and its conceptual content.

\begin{itemize}
	\item [a)] Umbral definition
		\begin{dmath}\label{eq:gaussian_cos}
			C_g(x) = {
				\frac{1}{1+\hat c x^2} \, \varphi_0,
				\quad
				\hat c^\alpha \varphi_0 = \frac{1}{\Gamma(1+ \alpha)}
			}
		\end{dmath}.
		In this formula, $\hat c$ denotes the umbral operator and $\varphi_0$ the corresponding vacuum.
		
	\item[b)] The Lorentzian umbral expression is mapped to the Gaussian. This is checked by expanding $C_g(x)$ in series and using the action of $\hat c$ on the vacuum:
		\begin{dmath*}
			C_g(x) = {
				\sum_{r=0}^\infty x^{2r}(-\hat c)^r \varphi_0 =
				\ee^{-x^2}
			}
		\end{dmath*}.
\end{itemize}

The umbral formalism therefore establishes a correspondence between rational functions and transcendental functions. Applying the partial fraction decomposition of the Lorentzian gives the identity
\begin{dmath}
	C_g(x) = \frac12 \left[\frac{1}{1 - i \hat c^{\frac12} x} + \frac{1}{1 + i \hat c^{\frac12} x} \right] \varphi_0
\end{dmath},
which suggests the introduction of the Gaussian exponential
\begin{dmath}
	e_g(x) = {
		\frac{1}{1 - i\hat c^{\frac12} x} \varphi_0 =
		\sum_{r=0}^\infty (i x)^{r}\hat c^{\frac{r}{2}} \, \varphi_0 =
		\sum_{r=0}^\infty \frac{(ix)^{r}}{\Gamma\left(\frac{r}{2} + 1 \right)}
	}
\end{dmath}
and of the associated Gaussian sine and cosine functions
\begin{dmath}
	S_g(x) = {
		\frac{1}{2i} \left[e_g(ix) - e_g(-ix)\right] \varphi_0 =
		\ee^{-x^2} \operatorname{erfi}(x)
	}
\end{dmath},
\begin{dmath}
	C_g(x) = {
		\frac{1}{2} \left[e_g(ix) + e_g(-ix)\right] \varphi_0 =
		\ee^{-x^2}
	}
\end{dmath}.
These functions are not trigonometric functions in the classical sense, since they do not possess the standard circular properties, such as periodicity or the ordinary derivative cycle of sine and cosine. Nevertheless, they share important structural analogies with the circular functions, and these analogies will be developed in the following sections.

The function $e_g(x)$ plays a central role in the partial fraction decomposition of the Lorentzian umbral image of the Gaussian. It is recognised as a particular case of the two-parameter Mittag--Leffler function \cite{WeissteinMittagLeffler}:
\begin{dmath}
	E_{\alpha,\,\beta}(x) = \sum_{r=0}^{\infty} \frac{x^r}{\Gamma(\alpha r + \beta)}
\end{dmath}.
For $\alpha = 1/2$ and $\beta = 1$, one has \cite{Haubold2011}
\begin{dmath}
	E_{\frac12,1}(x) = \ee^{x^2} \operatorname{erfc}(-x)
\end{dmath}.
The Mittag--Leffler function is often used in data analysis to interpolate between Lorentzian and Gaussian profiles, and it is a well-known generalisation of the exponential function. Moreover, the resemblance of the above construction with Euler's formula justifies, at least at the formal level, the terminology of Gaussian trigonometric functions.

The function $S_g(x)$, also related to the Faddeeva function \cite{vnF54}, plays an important role in plasma physics, where it is closely connected with the plasma dispersion function. The Gaussian sine and cosine also satisfy Kramers--Kronig-type relations:
\begin{dmath}
	S_g(x) = {
	-\frac{1}{\pi} \mathcal P \int_{-\infty}^\infty \frac{C_g(\xi)}{\xi - x}\,\mathrm d \xi,
	\quad
	C_g(x) = \frac{1}{\pi} \mathcal P \int_{-\infty}^\infty \frac{S_g(\xi)}{\xi - x}\,\mathrm d \xi
	}
\end{dmath}.
The role of the Gaussian exponential will be further discussed in the forthcoming sections, where we derive additional consequences of the formalism outlined above. In the final part of the paper, the same cyclic decompositions are reconsidered within the analytic umbral framework recently proposed in \cite{Riccianalytic2026} by one of the present authors (R.R.): the formal rational factorisations are realised as Mellin--Barnes pairings, thereby clarifying the origin of the component expansions, their contour representations, and their sectorial asymptotic behaviour.

\section{Cubic exponentials and Gaussian trigonometry}\label{sec:cubic_exp}

The present section remains within the formal umbral language. Its purpose is to show how the Gaussian case is only the quadratic member of a finite cyclic decomposition associated with roots of unity. The cubic case is treated explicitly because it is the first genuinely higher-order example and anticipates the analytic interpretation developed in \cref{sec:analytic_framework}.

In the previous section we have viewed at the trigonometric functions as the exponentiation of the square roots of the negative unit. In the following we refer at them as second order trigonometric functions. The third order partners are those obtained using the cubic roots of 1 and are defined as
\begin{dmath}
	C_0^3(x) = \frac13 \left[\ee^x + \ee^{\omega x} + \ee^{\omega^2 x}\right]
\end{dmath},
\begin{dmath}
	C_1^3(x) = \frac13 \left[\ee^x + \omega\ee^{\omega x} + \omega^2\ee^{\omega^2 x}\right]
\end{dmath},
\begin{dmath}
	C_2^3(x) = \frac13 \left[\ee^x + \omega^2\ee^{\omega x} + \omega\ee^{\omega^2 x}\right]
\end{dmath}.
where we have used $\omega = \ee^{i\frac{2\pi}{3}}$, $\,\omega^2 + \omega = -1$, $\,\omega^3 = 1$. It is accordingly evident that \cite{Babusci2012}
\begin{dmath}
	\frac{\mathrm d}{\mathrm d x} C_0^3(x) = {
		C_1^3(x),
		\quad 
		\frac{\mathrm d^2}{\mathrm d x^2} C_0^3(x) = C_2^3(x),
		\quad
		\frac{\mathrm d^3}{\mathrm d x^3} C_0^3(x) = C_0^3(x)
	}
\end{dmath}.

In the previous section we have touched different aspects of the Gaussian trigonometric functions, which has opened various elements of discussion. In our opinion the use of the factorisation of the Gaussian functions is a fairly natural consequence of the rational function we have chosen as image of the Gaussian function.

An obvious extension of the previous point of view is its application to the cubic exponential, namely
\begin{dmath}\label{eq:cubic_exp}
	\frac{1}{1+\hat c x^3} \, \varphi_0 = \ee^{-x^3}
\end{dmath},
where the umbral operator and the vacuum state have the same meaning as before.

The use of a standard procedure allows the following decomposition of the l.h.s. side of \cref{eq:cubic_exp}:
\begin{dmath}\label{eq:cubic_exp_decomposition}
	\frac{1}{1+\hat c x^3} \, \varphi_0 = \frac{1}{3}\sum_{r=0}^2 \frac{1}{1+\omega^r\hat c^{\frac13} x} \, \varphi_0
\end{dmath}
where $\omega$ is the same as before.

The conclusion we may draw from \cref{eq:cubic_exp,{eq:cubic_exp_decomposition}} is that
\begin{dmath}
	\ee^{-x^3} = {
		\frac{1}{3} \left[E_{\frac13,1}(-x) + E_{\frac13,1}(-\omega x) +E_{\frac13,1}(-\omega^2 x)\right] \varphi_0
	}
\end{dmath}.
The correctness of the last identity can be obtained by a direct numerical check.

The third order trigonometry is associated with the cubic roots of the unity \cite{Babusci2012}. We anticipate that there are three functions playing the role of the 

\begin{dmath}
	C_0^{g_3}(x) = \frac{1}{3} \left[E_{\frac13,1}(-x) + E_{\frac13,1}(-\omega x) +E_{\frac13,1}(-\omega^2 x)\right] \varphi_0
\end{dmath},
\begin{dmath}
	C_1^{g_3}(x) = \frac{1}{3} \left[E_{\frac13,1}(-x) + \omega^2 E_{\frac13,1}(-\omega x) + \omega E_{\frac13,1}(-\omega^2 x)\right] \varphi_0
\end{dmath},
\begin{dmath}
	C_2^{g_3}(x) = \frac{1}{3} \left[E_{\frac13,1}(-x) + \omega E_{\frac13,1}(-\omega x) + \omega^2 E_{\frac13,1}(-\omega^2 x)\right] \varphi_0
\end{dmath},
where $g_3$ denotes the correspondence with the cubic exponential. In this notation, we have the identifications $C_0^{g_2}(x) = C_g(x)$ and $C_1^{g_2}(x) = S_g(x)$.

The explicit expression of the functions in eq. (11) can be expressed as

\begin{dmath}
	C_0^{g_3}(x) = {E_{1,1}(-x^3) = \ee^{-x^3}}
\end{dmath},
\begin{dmath}
	C_1^{g_3}(x) = -x E_{1,\frac43}(-x^3)
\end{dmath},
\begin{dmath}
	C_2^{g_3}(x) = x^2 E_{1,\frac53}(-x^3)
\end{dmath}.
This result is a direct extension of what we already obtained for the Gaussian case. The last two identities have been derived by summing the series corresponding to the MLF in \cref{eq:cubic_exp}. The procedure is rather cumbersome and it will not be reported here because it implies boring technicalities. They have been verified both numerically and with the aid of Mathematica{\texttrademark}.

It is matter of a straightforward computation to infer that a super-Gaussian of order 4 can be decomposed using two different strategy
\begin{itemize}
	\item[a)] Decompose the quadratic parts only, namely
		\begin{dmath}
			\frac{1}{1+ax^4} =\frac12 \left(\frac{1}{1 - i a^{\frac12} x^2} + \frac{1}{1 + i a^{\frac12} x^2} \right)
		\end{dmath}
		Proceeding as before, we find
		\begin{dmath}
			C_g(x^2) = \frac12 \left[\frac{1}{1 - i \hat c^{\frac12} x^2} + \frac{1}{1 + i \hat c^{\frac12} x^2} \right] \varphi_0
		\end{dmath}
		\begin{dmath}
			S_g(x^2) = {
				\frac{1}{2i} \left[\frac{1}{1 - i \hat c^{\frac12} x^2} - \frac{1}{1 + i \hat c^{\frac12} x^2} \right] \varphi_0 =
				\ee^{-x^4} \operatorname{erfi}(x^2)
			}
		\end{dmath}

	\item[b)] Use the full decomposition
		\begin{dmath}
			\frac{1}{1+ax^4} =\frac14 \left(\frac{1}{1 - i a^{\frac14} x} + \frac{1}{1 + i a^{\frac14} x} + \frac{1}{1 - a^{\frac14} x} + \frac{1}{1 + a^{\frac14} x}\right)
		\end{dmath}
		and eventually end up with  
		\begin{dmath}
			\ee^{-x^4} = \frac12 \left(c^{g_4}(x) + ch^{g_4}(x) \right)
		\end{dmath},
		\begin{dmath}
			c^{g_4}(x) = {
				\frac12 \left(E_{\frac14,1}(ix) + E_{\frac14,1}(-ix)\right),
				\;
				ch^{g_4}(x) =
				\frac12 \left(E_{\frac14,1}(x) + E_{\frac14,1}(-x)\right)
			}
		\end{dmath}.
		\end{itemize}

The procedure we have followed mixes up roots of the negative and positive unit, a unifying formalism will be sketched in the final section.

We have so far used the umbral methodology to end up with the factorisation of exponentials and we have established a correspondence between a combination of Mittag--Leffler functions in the form of higher order trigonometric functions.

\section{Gaussian trigonometry and corresponding Fourier Transform}\label{sec:gaussian_trig}

This section records the formal Fourier-transform consequences of the cyclic umbral substitution. The calculations are intentionally kept at the operational level: the analytic restrictions on contour deformation, growth and admissible test functions are discussed later in \cref{sec:analytic_framework}.

According to the previous discussion the following correspondence can be envisaged
\begin{dmath}
	\ee^x \mapsto {E_{\frac12,1}(x) = e^{g_n}(x)}.
\end{dmath}
We can therefore ask whether such a correspondence can be exploited to state the existence of Mittag--Leffler based Fourier transforms:
\begin{dmath}\label{eq:ML_Fourier_transform}
	\mathcal F_{\frac1n}[f](k) = \int_{-\infty}^\infty f(x) \, E_{\frac1n, 1}(ikx) \, \mathrm d x
\end{dmath}.
We consider the case $n=2$ first and start from a classical example of ordinary Fourier transform

The procedure of evaluation of the transform in \cref{eq:ML_Fourier_transform} is exemplified below. We preliminary recall the ordinary Fourier transform of $f(x) = \ee^{-|x|}$ (see [12]):
\begin{dmath}
	\tilde f(k) = {
		\int_{-\infty}^\infty \ee^{-|x|} \, \ee^{-ikx} \, \mathrm d x =
		\frac{2}{1 + k^2}
	}
\end{dmath}.
Regarding the g-Fourier exponential we find
\begin{dmath}
	\mathcal F_{\frac12}[\ee^{-|x|}](k) = \left(\int_{-\infty}^\infty \ee^{-|x|} \frac{1}{1+ik\hat c^{\frac12} x} \, \mathrm d x\right) \varphi_0
\end{dmath}.
The explicit evaluation of the integral yields
\begin{dmath}
	\mathcal F_{\frac12}[\ee^{-|x|}](k) = {
		\left[\int_{-\infty}^\infty \ee^{-|x|} \left(\int_0^\infty \ee^{-s+isk\hat c^{\frac12} x}\,\mathrm d s\right) \, \mathrm d x\right] \varphi_0
	}
\end{dmath}
Interchanging the order of the integrals,  we find
\begin{dmath}
	\mathcal F_{\frac12}[\ee^{-|x|}](k) = {
		\left[2\int_{-\infty}^\infty \ee^{-s} \left(\int_0^\infty \ee^{-|x|+isk\hat c^{\frac12} x}\,\mathrm d x\right) \, \mathrm d s\right] \varphi_0
	}
\end{dmath}
Therefore, on account of \cref{eq:ML_Fourier_transform} we obtain
\begin{dmath}
	\mathcal F_{\frac12}[\ee^{-|x|}](k) = {
		\left(2\int_{-\infty}^\infty \ee^{-s} \frac{1}{1+\hat c (sk)^2}\, \mathrm d s \right) \varphi_0
		}
\end{dmath},
and finally, using \cref{eq:gaussian_cos}, we end up with
\begin{dmath}
	\mathcal F_{\frac12}[\ee^{-|x|}](k) = {
		2\int_{-\infty}^\infty \ee^{-s} \ee^{-(ks)^2}\, \mathrm d s
		}
\end{dmath}.
On the other hand, the same result is obtained after noting that $\ee^{-|x|}$ is an even function so that its g-Fourier transform can be written as
\begin{dmath}
	\int_{-\infty}^\infty \ee^{-|x|} \RE\left(\frac{1}{1+ik\hat c^{\frac12} x}\right)\varphi_0\,\mathrm d x =
	2\int_{-\infty}^\infty \ee^{-x}\ee^{-\omega x^2}\,\mathrm d x =
	\frac{\pi}{|k|}\ee^{\left(\frac{1}{2k}\right)^2}\operatorname{erfc}\left(\frac{1}{2|k|}\right)
\end{dmath}.

The definition of the g-Fourier transform, as it has been introduced, is incomplete and the concepts associated with the anti-transform counterpart should be fixed.

In the previous example we have used a result from the ordinary Fourier transform and used the umbral correspondence to state the corresponding $g_2$ transform a different strategy is outlined below, for the case
\begin{dmath}\label{eq:ML_Fourier_transform_gaussian}
	\mathcal F_{\frac12}[\ee^{-x^2}](k) = \int_{-\infty}^\infty \ee^{-x^2} \, E_{\frac12, 1}(ikx) \, \mathrm d x
\end{dmath}.
We introduce the umbral operator $\hat d$ acting on the vacuum as
\begin{dmath}
	\hat d^r = \frac{r!}{\Gamma\left(\frac{r}{2} + 1\right)}
\end{dmath}
We are therefore using an exponential as the umbral image of a Mittag--Leffler function, i.e. 
\begin{dmath}
	E_{\frac12,1}(ikx) = {
	\sum_{r=0}^\infty \frac{(ikx)^r}{\Gamma\left(\frac{r}{2} + 1\right)} =
	\sum_{r=0}^\infty \frac{(ikx)^r}{r!}\, \hat d^r \varphi_0 =
	\ee^{ikx\hat d} \varphi_0
	}
\end{dmath}.
The integral in \cref{eq:ML_Fourier_transform_gaussian} writes
\begin{dmath}
	\int_{-\infty}^\infty \ee^{-ax^2} \, E_{\frac12, 1}(ikx) \, \mathrm d x =
	\left(\int_{-\infty}^\infty \ee^{-ax^2}\ee^{ikx\hat d} \, \mathrm d x \right)\varphi_0
\end{dmath}.
The use of the rules of Gaussian integrals
\begin{dmath}
	\int_{-\infty}^\infty \ee^{-ax^2} \ee^{iqx} \, \mathrm d x =
	\sqrt{\frac{\pi}{a}}\ee^{-\frac{q^2}{4a}}
\end{dmath}
yields
\begin{dmath}
	\int_{-\infty}^\infty \ee^{-ax^2} \, E_{\frac12, 1}(ikx) \, \mathrm d x =
	\sqrt{\frac{\pi}{a}}\ee^{-k^2\hat{d}^2/4a}\varphi_0
\end{dmath}
After using the following identities,
\begin{dmath}
	\ee^{-k^2\hat{d}^2/4a}\varphi_0 = \sum_{r=0}^\infty \frac{1}{r!}\left(-\frac{k^2}{4a}\right) \hat d^{2r} \varphi_0
\end{dmath},
and
\begin{dmath}
	\hat d^{2r} \varphi_0 = \frac{(2r)!}{r!}
\end{dmath}.
we eventually find
\begin{dmath}
	\ee^{-k^2\hat{d}^2/4a}\varphi_0 = \sum_{r=0}^\infty \left(-\frac{k^2}{4a}\right)\frac{(2r)!}{(r!)^2}
\end{dmath}

Using the combinatorial identities
\begin{dmath}
	\frac{(2r)!}{(r!)^2} = {
		\binom{2r}{r},
		\quad
		\sum_{r=0}^\infty \binom{2r}{r} z^r =
		\frac{1}{\sqrt{1-4z}}
	}
\end{dmath}
the final result is
\begin{dmath}
	\int_{-\infty}^\infty \ee^{-ax^2} \, E_{\frac12, 1}(ikx) \, \mathrm d x =
	\sqrt{\frac{\pi}{a}}\frac{1}{\sqrt{1+k^2/a}}
\end{dmath}.

This result concludes our preliminary discussion on the Gaussian Fourier transform. In the following section, dedicated to final conclusions, we discuss some points requiring further comments and the plans for future works.

\section{Analytic framework}\label{sec:analytic_framework}

The preceding sections developed the cyclic decomposition in the formal language of indicial umbral calculus. We now reinterpret the same constructions within the analytic umbral framework introduced in \cite{Riccianalytic2026}. In this setting the umbral action is realised by a Mellin--Barnes pairing between a jump kernel and a ground state. This has two advantages. First, it explains why the formal exponential--rational transmutation produces Mittag--Leffler functions. Second, it shows that the cyclic components are not merely formal projections of a series, but analytic contour objects whose local and asymptotic expansions arise from different deformations of the same Mellin--Barnes integral. We remand to \cite{Riccianalytic2026} for more details.

The aim of this section is deliberately modest. We analyse explicitly the cases $n=2$ and $n=3$, because these are the cases used in the previous formal discussion, and then record the general $n$-fold pattern. The longer theory of asymptotics, Stokes phenomena and ramified companions is left to a separate treatment.

\subsection{Spectral transmutation and Mellin dilation}

The analytic counterpart of the formal identity used above is
\begin{dmath}\label{eq:identity}
	\left\langle
			\ee^{\eta\ee^u},
			\frac{\Gamma(1+t)}{\Gamma(1+t/n)}
		\right\rangle
	= {
	\left\langle
		\frac{1}{1-\eta\ee^{u/n}},
		\frac{1}{\Gamma(1+t)}
	\right\rangle = 
	E_{1/n,1}(\eta)
	}
\end{dmath}.
The first identity should be understood as the combination of two operations: a transmutation from the exponential Borel functional to the rational Borel functional, and a Mellin dilation of the spectral variable.

\noindent As proved in \cite{Riccianalytic2026} in the analytic framework and already evident at the formal level, given any admissible ground state function $\varphi(t)$, the following identity holds:
\begin{dmath*}
	\left\langle
			\ee^{\eta\ee^u},
			\varphi(t)
		\right\rangle
	= {
	\left\langle
		\frac{1}{1-\eta\ee^{u}},
		\frac{\varphi(t)}{\Gamma(1+t)}
	\right\rangle 
	}
\end{dmath*}.
We say that the pairing of the exponential Borel functional with $\varphi(t)$ is equivalent to the pairing of the rational Borel functional with the \emph{transmuted} ground state $\varphi(t)/\Gamma(1+t)$. This is an immediate consequence of a spectral duality holding between the Mellin representatives (a.k.a. \emph{jump kernels}, see below) of the two Borel functionals, namely:
\begin{dmath*}
	J_{\mathrm{exp}}(t; \zeta) = \frac{J_{\mathrm{rat}}(t; \zeta)}{\Gamma(1+t)}	
\end{dmath*}.
\Cref{eq:identity} is obtained by combining this transmutation with a further dilation of the Mellin--Barnes variable $t$, as clarified below.

With the Mellin convention
\begin{dmath*}
	\mathcal M[f](s)
	=
	\int_0^\infty z^{s-1}f(z)\,\mathrm d z
\end{dmath*},
the jump kernels of the exponential and rational Borel functionals are
\begin{dmath*}
	J_{\mathrm{exp}}(t;\eta) = {
		\mathcal M_z[\ee^{\eta z}](-t)
		=
		(-\eta)^t\Gamma(-t)
	}
\end{dmath*}
and
\begin{dmath*}
	J_{\mathrm{rat}}(t;\eta) = {
		\mathcal M_z\left[\frac{1}{1-\eta z}\right](-t)
		=
		-(-\eta)^t\pi\csc(\pi t)
	}
\end{dmath*},
with the branch fixed by the chosen sector. Therefore the analytic pairing corresponding to the formal cyclic exponential is
\begin{dmath*}
	I_n(\eta)
	= {
		\left\langle
			\ee^{\eta\ee^u},
			\frac{\Gamma(1+t)}{\Gamma(1+t/n)}
		\right\rangle
		=
		\frac{1}{2\pi i}
		\int_{\mathcal C_t}
		(-\eta)^t\Gamma(-t)
		\frac{\Gamma(1+t)}{\Gamma(1+t/n)}
		\,\mathrm d t
	}
\end{dmath*}.
Using $\Gamma(-t)\Gamma(1+t)=-\pi/\sin(\pi t)$, this becomes
\begin{dmath*}
	I_n(\eta)
	=
	-
	\frac{1}{2\pi i}
	\int_{\mathcal C_t}
	(-\eta)^t
	\frac{\pi}{\sin(\pi t)}
	\frac{1}{\Gamma(1+t/n)}
	\,\mathrm d t
\end{dmath*}.
Closing the contour to the right gives the primary expansion
\begin{dmath*}
	I_n(\eta)
	= {
		\sum_{m=0}^{\infty}
		\frac{\eta^m}{\Gamma(1+m/n)}
		=
		E_{1/n,1}(\eta)
	}
\end{dmath*}.
Thus the formal umbral object of the previous sections is analytically the Mittag--Leffler function $E_{1/n,1}$.

Performing the change of variable $t=ns$ yields
\begin{dmath*}
	I_n(\eta)
	=
	-
	\frac{n}{2\pi i}
	\int_{\mathcal C_s}
	(-\eta)^{ns}
	\frac{\pi}{\sin(n\pi s)}
	\frac{1}{\Gamma(1+s)}
	\,\mathrm d s
\end{dmath*}.
This is the Mellin--Barnes form of the rational pairing
\begin{dmath*}
	I_n(\eta)
	=
	\left\langle
		\frac{1}{1-\eta\ee^{u/n}},
		\frac{1}{\Gamma(1+t)}
	\right\rangle
\end{dmath*}.
This identity is the analytic version of the formal transformations used in the previous sections.

\subsection{The quadratic case}

For $n=2$ one obtains
\begin{dmath*}
	I_2(\eta)
	=
	E_{1/2,1}(\eta)
\end{dmath*},
with the closed form
\begin{dmath*}
	E_{1/2,1}(\eta)
	=
	\ee^{\eta^2}\operatorname{erfc}(-\eta)
\end{dmath*}.
The rational decomposition
\begin{dmath*}
	\frac{1}{1-\eta\ee^{u/2}}
	=
	\frac{1}{1-\eta^2\ee^u}
	+
	\frac{\eta\ee^{u/2}}{1-\eta^2\ee^u}
\end{dmath*}
induces the splitting
\begin{dmath*}
	I_2(\eta)
	=
	I_{2,0}(\eta)+I_{2,1}(\eta)
\end{dmath*},
where
\begin{dmath*}
	I_{2,0}(\eta)
	= {
		\sum_{q=0}^{\infty}
		\frac{\eta^{2q}}{\Gamma(1+q)}
		=
		\ee^{\eta^2}
	}
\end{dmath*}
and
\begin{dmath*}
	I_{2,1}(\eta)
	= {
		\sum_{q=0}^{\infty}
		\frac{\eta^{2q+1}}{\Gamma(q+3/2)}
		=
		\eta E_{1,3/2}(\eta^2)
	}
\end{dmath*}.
For the Gaussian choice $\eta=ix$, the even component gives the Gaussian factor
\begin{dmath*}
	I_{2,0}(ix)
	=
	\ee^{-x^2}
\end{dmath*},
whereas the odd component gives the Gaussian sine/Faddeeva companion
\begin{dmath*}
	\frac{1}{i}I_{2,1}(ix)
	=
	\ee^{-x^2}\operatorname{erfi}(x)
\end{dmath*}.
This is precisely the analytic counterpart of the Gaussian trigonometric functions introduced formally in \cref{sec:introduction}.

The same splitting is obtained from the Mellin--Barnes kernel through
\begin{dmath*}
	\frac{\pi}{\sin(2\pi s)}
	=
	\frac{1}{2}
	\left(
		\frac{\pi}{\sin(\pi s)}
		-
		\frac{\pi}{\sin\!\left(\pi\left(s-\frac{1}{2}\right)\right)}
	\right)
\end{dmath*}.
Thus the partial fraction factorisation of the rational umbral kernel and the spectral splitting of the Mellin kernel are the same operation seen in two languages.

\subsection{The cubic case}

For $n=3$ one has
\begin{dmath*}
	I_3(\eta)
	=
	E_{1/3,1}(\eta)
\end{dmath*}.
The rational decomposition is
\begin{dmath*}
	\frac{1}{1-\eta\ee^{u/3}}
	=
	\frac{1}{1-\eta^3\ee^u}
	+
	\frac{\eta\ee^{u/3}}{1-\eta^3\ee^u}
	+
	\frac{\eta^2\ee^{2u/3}}{1-\eta^3\ee^u}
\end{dmath*}.
Accordingly,
\begin{dmath*}
	I_3(\eta)
	=
	I_{3,0}(\eta)+I_{3,1}(\eta)+I_{3,2}(\eta)
\end{dmath*},
with
\begin{dmath*}
	I_{3,0}(\eta)
	= {
		\sum_{q=0}^{\infty}
		\frac{\eta^{3q}}{\Gamma(1+q)}
		=
		\ee^{\eta^3}
	}
\end{dmath*},
\begin{dmath*}
	I_{3,1}(\eta)
	= {
		\sum_{q=0}^{\infty}
		\frac{\eta^{3q+1}}{\Gamma(q+4/3)}
		=
		\eta E_{1,4/3}(\eta^3)
	}
\end{dmath*},
and
\begin{dmath*}
	I_{3,2}(\eta)
	= {
		\sum_{q=0}^{\infty}
		\frac{\eta^{3q+2}}{\Gamma(q+5/3)}
		=
		\eta^2E_{1,5/3}(\eta^3)
	}
\end{dmath*}.
Taking $\eta=-x$ gives the cubic Gaussian components displayed formally in Section 2:
\begin{dmath*}
	I_{3,0}(-x)
	=
	\ee^{-x^3}
\end{dmath*},
\begin{dmath*}
	I_{3,1}(-x)
	=
	-xE_{1,4/3}(-x^3)
\end{dmath*},
and
\begin{dmath*}
	I_{3,2}(-x)
	=
	x^2E_{1,5/3}(-x^3)
\end{dmath*}.
The Mellin--Barnes splitting is governed by
\begin{dmath*}
	\frac{\pi}{\sin(3\pi s)}
	=
	\frac{1}{3}
	\left(
		\frac{\pi}{\sin(\pi s)}
		-
		\frac{\pi}{\sin\!\left(\pi\left(s-\frac{1}{3}\right)\right)}
		+
		\frac{\pi}{\sin\!\left(\pi\left(s-\frac{2}{3}\right)\right)}
	\right)
\end{dmath*}.
This is the analytic origin of the three cyclic components obtained formally by the cubic roots of unity.

\subsection{General cyclic decomposition and Fourier-transform synthesis}

The preceding cases are the first members of the general decomposition
\begin{dmath*}
	\frac{1}{1-\eta\ee^{u/n}}
	=
	\sum_{r=0}^{n-1}
	\frac{\eta^r\ee^{ru/n}}{1-\eta^n\ee^u}
\end{dmath*}.
Consequently,
\begin{dmath*}
	I_n(\eta)
	=
	\sum_{r=0}^{n-1}
	I_{n,r}(\eta)
\end{dmath*},
where
\begin{dmath*}
	I_{n,r}(\eta)
	= {
		\eta^r
		\sum_{q=0}^{\infty}
		\frac{\eta^{nq}}{\Gamma\!\left(1+q+\frac{r}{n}\right)}
		=
		\eta^rE_{1,1+r/n}(\eta^n)
	}
\end{dmath*}.
Equivalently, if $\omega=\ee^{2\pi i/n}$, then
\begin{dmath*}
	I_{n,r}(\eta)
	=
	\frac{1}{n}
	\sum_{k=0}^{n-1}
	\omega^{-rk}I_n(\omega^k\eta)
\end{dmath*}.
This root-of-unity projector selects the powers $m\equiv r\pmod n$ in the primary expansion of $E_{1/n,1}(\eta)$.

The Fourier-transform construction of \cref{sec:gaussian_trig} can be summarised in the same analytic language. The primary definition is spectral:
\begin{dmath*}
	\mathfrak F_n[f](k)
	=
	\left\langle
		\mathcal F[f](k\ee^u),
		\frac{\Gamma(1+t)}{\Gamma(1+t/n)}
	\right\rangle
\end{dmath*}.
When the pairing may be interchanged with the spatial integral, this gives the kernel representation
\begin{dmath*}
	\mathfrak F_n[f](k)
	=
	\int_{-\infty}^{\infty}
	f(x)E_{1/n,1}(-ikx)
	\,\mathrm d x
\end{dmath*}.
The Gaussian example in \cref{sec:gaussian_trig} is a concrete instance of this principle. For $n=2$ and $f(x)=\ee^{-a x^2}$ one obtains
\begin{dmath*}
	\int_{-\infty}^{\infty}
	\ee^{-a x^2}E_{1/2,1}(ikx)
	\,\mathrm d x
	=
	\sqrt{\frac{\pi}{a+k^2}}
\end{dmath*},
with the branch fixed by analytic continuation from the disk of convergence of the corresponding power series.

The analytic framework therefore reorganises the formal construction as follows: the finite cyclic decompositions of the previous sections are root-of-unity splittings of a Mellin--Barnes kernel, while the associated Fourier transforms are spectral deformations of the ordinary Fourier transform by the same Gamma-ratio ground state.
\subsection{Analytic synthesis of the cyclic construction}\label{subsec:analytic_synthesis}

The analysis carries out in this section gives a contour-level interpretation of the formal cyclic decompositions developed in the preceding sections. The fundamental analytic object is the Mellin--Barnes pairing, reproducing and generalising the formal construction. We prove, in particular, that the quadratic and cubic examples are special cases of a general spectral transmutation between exponential and rational umbral functionals.

The cyclic components are analytic objects in the same sense. They may be obtained either from the root-of-unity factorisation of the rational kernel or from the splitting of the Mellin factor \(\pi/\sin(n\pi s)\). This gives the contour realisation of the higher-order trigonometric components introduced formally in \cref{sec:cubic_exp,sec:gaussian_trig}.

The same analytic point of view also clarifies the Fourier-transform construction, whose natural definition is spectral:
\begin{dmath*}
	\mathfrak F_n[f](k)
	=
	\left\langle
		\mathcal F[f](k\ee^u),
		\frac{\Gamma(1+t)}{\Gamma(1+t/n)}
	\right\rangle
\end{dmath*}.
Whenever the pairing can be exchanged with the spatial integral, this becomes a Fourier transform with Mittag--Leffler kernel. The distinction is important: the spatial representation is useful, but it may impose stronger convergence restrictions than the primary spectral pairing.

In summary, the analytic framework does not replace the formal umbral construction. Rather, it explains why the same cyclic functions arise from rational factorisation, root-of-unity projection, Mellin--Barnes splitting, and spectral deformation of Fourier transforms.

Several questions remain open. A full inversion theory for the umbral Fourier transform has not been developed here. The admissible classes of test functions and distributions should be characterised in the analytic-functional setting. Finally, the asymptotic and Stokes behaviour of the cyclic components, already visible from the Mellin--Barnes representation, deserves a separate investigation. These directions will be addressed in future work.

\section{Concluding Remarks and Perspectives}\label{sec:conclusion}

The analysis developed in the previous sections shows that the umbral formalism provides a consistent and operationally transparent framework for extending the notion of Fourier transform to a broader class of functions generated by non-standard exponentials. In particular, the introduction of the umbral operator, acting on the vacuum according to
\begin{dmath}
	\hat c^r \varphi_0 = \frac{1}{r + 1}
\end{dmath}
allows one to translate generalized (g-)exponentials into ordinary exponential structures, preserving the algebraic backbone of Fourier analysis.

The main result of this work is that the g-Fourier transform can be evaluated by mapping it onto a standard Fourier transform, provided that the umbral substitution is implemented with the correct vacuum prescription. As shown in the derivation leading to Eqs. (10)--(15), the consistency of the procedure relies on replacing factorial structures with their gamma-function counterparts, ensuring that the operational correspondence remains exact and yields the correct final expressions.

A natural and significant extension of the present framework is obtained by considering the general case of $\ee^{-x^n}$exponentials of the  which can be treated within the same umbral scheme. In this context, the decomposition in terms of the $n$-th roots of unity emerges as the unifying principle. Indeed, the factorisation of the corresponding algebraic structures allows one to express the generalised exponential as a superposition of elementary contributions associated with the roots $\omega_k = \ee^{2\pi i k/n}, \; k = 0, 1, \dots\ n-1$.

This leads to a natural splitting of the transform into components, each governed by a distinct branch determined by $\omega_k$. The previously discussed quadratic and quartic cases therefore appear as particular realisations of a general construction, where the umbral operator effectively encodes the combinatorial structure associated with the roots of unity. This observation suggests that the entire hierarchy of generalised Fourier transforms can be organised within a single algebraic framework based on cyclic symmetry.

Another point deserving attention concerns the inverse (anti-transform) procedure. While the forward g-Fourier transform is well defined through the umbral correspondence, the construction of its inverse requires additional care. In particular, one should clarify the role of the dual vacuum and the corresponding inverse operator, ensuring that the inversion reproduces the original function without ambiguity. It is plausible that the inverse transform can be formulated by exploiting the same roots-of-unity decomposition, leading to a reconstruction formula involving conjugate branches. However, a rigorous formulation of this procedure remains an open problem and should be addressed in future investigations.

A further perspective is offered by the possibility of extending the present formalism beyond polynomial exponentials. The structure uncovered here, especially in the $n$-th order case, suggests a connection with more sophisticated function classes, such as those related to elliptic functions. Since elliptic functions can be viewed as arising from doubly periodic structures, one may speculate that an appropriate generalisation of the umbral exponential---possibly involving multiple operators or higher-dimensional vacuum states---could accommodate such periodicities. In this sense, the decomposition in terms of roots of unity may represent a first step toward a more general theory where algebraic and periodic structures are simultaneously encoded.

Another possible extension of the formalism could be realised by applying the same cyclic decomposition, analytically interpreted as a spectral transmutation followed by a Mellin dilation, to other families of functions. Promising results have already been obtained in the case of Bessel and Bessel--Clifford functions and will be presented in a following article.

In summary, the umbral decomposition of Fourier transforms preserves the operational simplicity of the classical Gaussian case while extending it to a wide class of generalised exponentials. The emergence of the roots-of-unity structure provides a unifying viewpoint, suggesting that different cases can be treated within a single algebraic scheme. The formulation of the inverse transform, as well as the extension toward elliptic and more general functions, remain open and promising directions for future research.

\bibliographystyle{unsrt}
\bibliography{bibliography_UMFF}


\end{document}